\documentclass{gtart}
\gtart
\usepackage{amsthm}
\usepackage{amsgen}
\usepackage{amsmath,amssymb}
\usepackage{curvesls}
\usepackage{epic}
\usepackage[all]{xy}
\usepackage[dvips]{graphicx}
\usepackage{psfrag}
\usepackage{color}
\usepackage{ccaption}
\usepackage{multicol}

\newtheorem{thm}{Theorem}[section]

\newtheorem{prop}[thm]{Proposition}
\newtheorem{eg}[thm]{Example}

\theoremstyle{definition}

\newcommand{\Real}{{\mathbb R}}

\newcommand{\Rat}{{\mathbb Q }}
\newcommand{\Zed}{{\mathbb Z}}

\newcommand{\GL}{{\mathrm{GL}}}
\newcommand{\gO}{{\mathrm{O}}}
\newcommand{\gSO}{{\mathrm{SO}}}
\newcommand{\Sym}{{\mathrm{Sym}}}
\newcommand{\Gr}{{\mathrm{Gr}}}

\newcommand{\Sec}{{\Gamma^{\mathbb S}_1 VN}}

\newcommand{\Emb}{{\mathrm{Emb}}}

\newcommand{\tspinc}{{{spin}\raisebox{1mm}{\footnotesize c}}}
\newcommand{\tSpinc}{{{Spin}\raisebox{1mm}{\footnotesize c}}}
\newcommand{\mSpinc}{{\text{Spin}\raisebox{1mm}{\footnotesize c}}}
\newcommand{\mSpin}{{\text{Spin}}}

\begin{document}

\title{Combinatorial spin structures \\ on triangulated manifolds}

\author{Ryan Budney}
\address{Mathematics and Statistics, University of Victoria, \\
PO BOX 1700 STN CSC, Victoria, B.C., Canada V8W 2Y2}
\email{rybu@uvic.ca }

\begin{abstract}
This paper gives a combinatorial description of spin and \tspinc-structures
on triangulated manifolds of arbitrary dimension.  These encodings of 
spin and \tspinc-structures are established primarily for the purpose of aiding in 
computations. The novelty of the approach is that we rely heavily on the naturality of 
binary symmetric groups to avoid lengthy explicit constructions of smoothings
of PL manifolds. 
\end{abstract}

\primaryclass{57R15}
\secondaryclass{57R05, 55S35}
 
\keywords{spin structure, \tspinc \ structure, triangulation, manifold}

\maketitle

\section{Introduction}\label{INTRODUCTION}

In this paper a framework for combinatorially representing 
 spin and \tspinc-structures on triangulated manifolds in a manner suitable for 
computer implementation is built.  This should be seen as part of a general effort to 
merge the techniques of algorithmic 3-manifold theory, such as triangulations, 
normal surface theory and geometrization with elements of 4-manifold theory where
gauge-theoretic invariants often require additional structures. 

The governing perspective on spin and \tspinc-structures in this paper comes from 
 the obstruction-theoretic approach to spin structures of Milnor \cite{Milnor}.  Although 
 Milnor's approach is fundamentally combinatorial in nature, there is some non-trivial 
 work to translate Milnor's language into a language a modern computer can use.  To this 
 end, we put combinatorial spin structures in a formalism perhaps most comparable 
 to Forman's discrete Morse theory \cite{F}. It is assumed the reader is familiar 
 with obstruction theory on manifolds along the lines of Milnor and Stasheff \cite{MS}.  
 Other references like Whitehead \cite{Whitehead} and Gompf and Stipicz \cite{GS} are 
 also excellent resources for basic obstruction theory. 

Relatively flexible triangulations are allowed in this article.  For example, 
 {\it unordered delta complexes} \cite{Hat} suffice.  The {\it ideal triangulations} of 
 Thurston \cite{Thu1}, a further weakening of unordered delta complexes, are also 
 perfectly acceptable. Ideal triangulations are unordered delta complexes, such that if one
removes a finite collection of vertices, one obtains a manifold.  In short, a 
{\it triangulation} in this paper is a space constructed by gluing simplices together 
via affine-linear identifications of their boundary facets, and where we demand that the 
characteristic maps of every simplex is an embedding when restricted to the interior 
of the simplex. 

Readers comfortable with the basics of triangulations, spin structures and obstruction theory can jump to
Section \ref{mainsec} for the primary constructions of this paper. 
In the literature, there are several available tools for combinatorially representing $3$ and 
$4$-manifolds with additional structure on their tangent bundles.   The Kaplan Algorithm 
\cite{Kaplan} was perhaps the first (see \cite{GS} \S 5.6, 5.7 for a modern exposition). Kaplan's
 Algorithm gives a simple framework to represent spin structures on a $3$-manifold given by an
 integral surgery presentation, and provides a simple tool to determine when such spin-structures 
extend over the bounding $4$-manifold.  Another combinatorial representation of 
$3$-manifolds are {\it spines}, popularized by Matveev \cite{Mat}. Techniques to represent 
spin-structures on $3$-manifold and $4$-manifold spines were developed by Benedetti and 
Petronio \cite{BP, BP01, BP02, BP2}.  The techniques in this paper would 
be described as being in the language of the `frame along the dual $1$-skeleton' in 
\cite{BP2}. \tSpinc-structures on simplicially-triangulated $3$-manifolds can be described as 
the combinatorial Euler structures of Turaev \cite{Turaev}.  \'Etienne 
Gallais has recently used this technique to study combinatorial Euler structures on 
triangulated $3$-manifolds \cite{Gallais} 
using Forman's combinatorial vector fields to represent Euler structures.  One of Gallais's 
observations is that with these techniques, not all combinatorial Euler structures are 
represented on delta complexes.  Simplicial triangulations are required 
to capture {\it all} \tspinc-structures using this technique.  We wish to avoid simplicial 
triangulations, as unordered delta complexes have shown themselves to be rather
 efficient means for describing interesting manifold types in both 
3-manifold theory \cite{Thu1, BBP} and 4-manifold theory \cite{BBH, BH}.

\section{Notation, obstruction theory}\label{defsnot}

Throughout this paper, $N$ will be a PL $n$-manifold that will be endowed with
a triangulation or a CW-structure, often both.  If the cell structure is unambiguous,
the $i$-skeleton will be denoted by $N^i$.  

Given a fibre bundle $\psi : E \to B$ with fibre $F$, and a subspace $X \subset B$, 
the {\it restriction bundle} is the map $\psi_{|\psi^{-1}(X)} : \psi^{-1}(X) \to X$ 
which also has fibre $F$.  
We abbreviate $\psi_{|\psi^{-1}(X)}$ with $\psi_{|X}$. 

A {\it trivialization} of a vector bundle $\psi : E \to B$ is an ordered $k$-tuple of
vector fields that form a basis for each and every fibre.  Trivializations 
correspond to vector bundle isomorphisms $B \times \Real^k \to E$ via the map
$(b,x_1,\cdots,x_k) \longmapsto \sum_{i=1}^k x_i \vec v_i(b)$ where 
$(\vec v_i : B \to E)_{i \in \{ 1,2,\cdots,k\}}$ is the trivialiation.

A vector bundle $\psi : E \to N$ is {\it orientable} if and only if there is a 
trivialization of $\psi_{|N^1}$. Given a trivialization of $\psi_{|N^1}$, the homotopy class of its 
restriction to $N^0$ is called an {\it orientation} of $\psi$. If a vector bundle $\psi : E \to N$ is 
orientable, its set of orientations admits a free transitive action of $H^0(N,\Zed_2)$ -- the action is 
given by flipping orientations on path-components of $N$.  

In the language of classifying maps a vector bundle $\psi : E \to N$ is {\it orientable} if and only if its 
classifying map $N \to \Gr_{\infty,k} \equiv B\gO_k$ lifts to the Grassmannian of oriented $k$-subspaces of 
$\Real^\infty$, $\Gr^+_{\infty,k} \equiv B\gSO_k$. 

$$\xymatrix{ 
         & B\gSO_k \ar[d] \\
N \ar[ur] \ar[r] & B\gO_k } $$

An {\it orientation} of $N$ is the homotopy class of this lift.  The fact that this is equivalent 
to the previous definition is described in the references \cite{Milnor, MS}. The key ingredient in 
this interpretation is that $\gSO_k$ is the path-component of the identity in $\gO_k$.  This implies 
that the exact CW-structure on the space $N$ is not relevant to the 
existence of orientations, which is one reason to prefer this formalism.  If $N$ is a smooth manifold, 
{\it orientability and orientations of $N$} refer to orientability and orientations 
of the tangent bundle $\pi : TN \to N$. 

The {\it $n^\text{th}$ spin group} we denote by $\mSpin_n$.  This is defined as is the unique connected Lie group 
which admits an onto $2:1$ Lie group homomorphism $\mSpin_n \to \gSO_n$. Since $\pi_1 \gSO_n$ is 
cyclic of order $2$ or infinite-cyclic, this is well-defined.  A vector bundle $\psi : E \to N$ 
admits a spin-structure if the classifying map $N \to B\gO_k$ admits a lift $N \to B\mSpin_k$.  
A spin structure is a homotopy class of map $N \to B\mSpin_k$ such that the composite with 
$B\mSpin_k \to B\gO_k$ is a classifying map for the bundle $\psi$.  Since the homomorphism 
$\mSpin_k \to \gO_k$ factors as a composition $\mSpin_k \to \gSO_k \to \gO_k$, spin structures 
induce orientations. 

$$\xymatrix{ & B\mSpin_k \ar[d] \\
N \ar[ur] \ar[r] \ar[dr] & B\gSO_k \ar[d] \\
             & B\gO_k }$$

Since $\pi_1 \gSO_k \simeq \Zed_2$ for $k \geq 3$, the corresponding description for spin structures in 
the obstruction-theoretic setting is that $\psi : E \to N$ {\it admits a spin structure} if and only 
if there exists a trivialization of $\psi_{|N^2}$.  Given such a trivialization, the homotopy class 
of its restriction to $N^1$ is a {\it spin structure}.  The case $k=2$ is special since $\pi_1 \gSO_2$ 
is infinite cyclic. Typically in the literature people phrase the obstruction-theoretic formulation as 
saying $\psi \oplus \epsilon^1$ admits a spin structure, where $\epsilon^1 : N \times \Real \to N$ is 
the trivial $1$-dimensional bundle over $N$, but one could just as easily describe it in terms of 
trivializations of $\psi_{|N^1}$ such that the obstructions to extending over $N^2$ are all even. 

The $k^\text{th}$ {\it complex spin group}, $\mSpinc_k$ is the group 
$(\mSpin_k \times \mSpin_2)/\Zed_2 \equiv \mSpin_k \times_{\Zed_2} \mSpin_2$. This means we are taking the product 
of the $k^\text{th}$ spin group with the $2^\text{nd}$ spin group, and modding out by one copy of $\Zed_2$ acting 
diagonally on the product via the covering action on the respective spin groups.  Via projection to the right
and left factor respectively 
this group admits two extensions: $\mSpin_k \to \mSpinc_k \to \gSO_2$ and $\mSpin_2 \to \mSpinc_k \to \gSO_k \equiv \mSpin_k/\Zed_2$.  
The latter extension is used to define \mSpinc-structures, and the former gives the inclusion $\mSpin_k \to \mSpinc_k$.  

$$\xymatrix{ & B\mSpin_k \ar[d] \\
             & B\mSpinc_k \ar[d] \\
N \ar[uur] \ar[ur] \ar[r] \ar[dr] & B\gSO_k \ar[d] \\
             & B\gO_k }$$

A vector bundle $\psi : E \to N$ admits a \tspinc-structure if the classifying map $N \to B\gO_k$ admits 
a lift to $B\mSpinc_k$ \cite{gom}.  A {\it \tspinc-structure} is a homotopy class of map 
$N \to B\mSpinc_k$ such that the composition with $B\mSpinc_k \to \gO_k$ classifies the bundle $\psi$.  
To interpret a \tspinc-structure, notice that if one composes with the former extension, one gets a 
map $N \to B\gSO_2$ which classifies an oriented $2$-dimensional vector bundle over $N$. Alternatively 
this is a $1$-dimensional $\mathbb C$-bundle over $N$.  If $\nu : E' \to N$ is the $1$-dimensional 
$\mathbb C$-bundle over $N$ classified by this map, then $\psi \oplus \nu : E \oplus E' \to N$ is 
classified by the corresponding map $N \to B\gSO_k \times B\gSO_2 \equiv B(\gSO_k \times \gSO_2)$.  
Consider $\gSO_k \times \gSO_2$ as a subgroup of $\gSO_{k+2}$.  This group is covered by some subgroup of 
$\mSpin_{k+2}$, and by design this group is isomorphic to $\mSpin_k \times_{\Zed_2} \mSpin_2$.  Thus 
a \mSpinc-structure on a bundle $\psi : E \to N$ consists of two things: a complex line 
bundle $\nu : E' \to N$ and a spin-structure on $\psi \oplus \nu$.   Given this, \tspinc-structures 
can be readily transcribed into an obstruction-theoretic formalism.  A complex line bundle is 
classified by a map $N \to B\gSO_2 \equiv K(\Zed,2)$ and homotopy classes of maps $N \to K(\Zed,2)$ are 
in bijective correspondence with elements of $H^2(N, \Zed)$.  Thus a \tspinc-structure on $N$ is 
prescribed by such a cohomology class, together with a homotopy class of trivialization of 
$(\psi \oplus \nu)_{|N^1}$ which extends to $N^2$. 

When working with a triangulation $T$ of a manifold $N$, we will make heavy usage of the {\it dual polyhedral 
decomposition}.  This construction originated in the work of Poincar\'e, and is available in \cite{ST}. Since 
these ideas are no longer in wide circulation and we need some fixed notation to refer to this 
decomposition, a brief sketch is given.  Denote the standard $n$-simplex by 
 $$\Delta^n = \{(x_0,\cdots,x_n) \in \Real^{n+1} : x_i \geq 0 \ \forall i \text{ and } x_0+x_1+\cdots+x_n = 1\}.$$
  
For $i \in \{0,1,\cdots,n\}$ the $i$-th face map of $\Delta^n$ is $f_i : \Delta^{n-1} \to \Delta^n$ given by 
$f_i(x_0,\cdots,x_{n-1}) = (x_0,x_1,\cdots,x_{i-1}, 0, x_i, x_{i+1}, \cdots, x_{n-1})$. 
Given a permutation $\sigma \in \Sigma_{n+1} \equiv \Sigma(\{0,1,\cdots,n\})$, the induced automorphism of $\Delta^n$ 
is denoted $\sigma_* : \Delta^n \to \Delta^n$ and is defined by
$\sigma_*(x_0,x_1,\cdots,x_n) = (x_{\sigma^{-1}(0)}, x_{\sigma^{-1}(1)}, \cdots, x_{\sigma^{-1}(n)})$.  
An {\it unordered delta complex} is a CW-complex $X$ such that the domains of the attaching maps are the boundaries of 
simplices (rather than discs), $\phi : \partial \Delta^n \to X^{(n-1)}$, and for each $i$, the composite 
satisfies $\phi \circ f_i = \Phi \circ \sigma_*$ where $\Phi : \Delta^{n-1} \to X^{(n-1)}$ is the 
characteristic map of some $(n-1)$-simplex, and $\sigma \in \Sigma_n$ is 
some permutation. If all the permutations $\sigma$ were the identity, $X$ would be an {\it ordered delta complex}. 

Let $[0,n] = \{0,1,\cdots,n\}$, and let $I$ denote a subset of $[0,n]$. The {\it dual polyhedral bit} $\delta_I$ 
of $\Delta^n$ is the convex hull of the barycentres of all faces of $\Delta^n$ with vertex-sets a super-set of $I$.
Thus, $\delta_{[0,n]}$ is the barycentre of $\Delta^n$ and $\delta_{[0,n] \setminus \{i\}}$ is the convex hull of the barycentre of $\Delta^n$ together with the barycentre of the $i$-th face of $\Delta^n$.  One can define $\delta_I$ via a system of equations, as well

$$\delta_I = \{(x_0,x_1, \cdots, x_n) \in \Delta^n : x_i \geq x_j \ \forall i \in I \text{ and } j \in [0,n] \}.$$

If $T$ is a triangulation of a manifold $N$, and $\chi : \Delta^n \to N$ the characteristic map of a simplex, 
$\chi(\delta_I)$ is defined to be a {\it dual polyhedral bit of the triangulation $T$}. Given an $i$-dimensional 
simplex $\sigma$ of $T$, {\it the closed dual $(n-i)$-cell} corresponding to $\sigma$ is the union 
of all $(n-i)$-dimensional dual polyhedral bits corresponding to $\sigma$ in all the top-dimensional simplices
containing $\sigma$. The collection 
of all dual cells forms a CW-decomposition of $N$, called {\it the polyhedral decomposition of $N$ dual to $T$}.  
We denote this dual polyhedral decomposition by $P$ throughout the paper.   Given a triangulation $T$ or
CW-complex $P$, we denote the set of $k$-cells by $T_k$ and $P_k$ respectively, while the $k$-skeleton
we continue to denote by $T^k$ and $P^k$ respectively.  The key feature of the dual decomposition is that
for every $i$-simplex $\sigma \in T_i$ there is one and only one dual $(n-i)$-cell $e^{n-i} \in P_{n-i}$ with 
$\sigma \cap e^{n-i} \neq \emptyset$. The non-empty intersection is the barycentre of $\sigma$. 

$$\includegraphics[width=8cm]{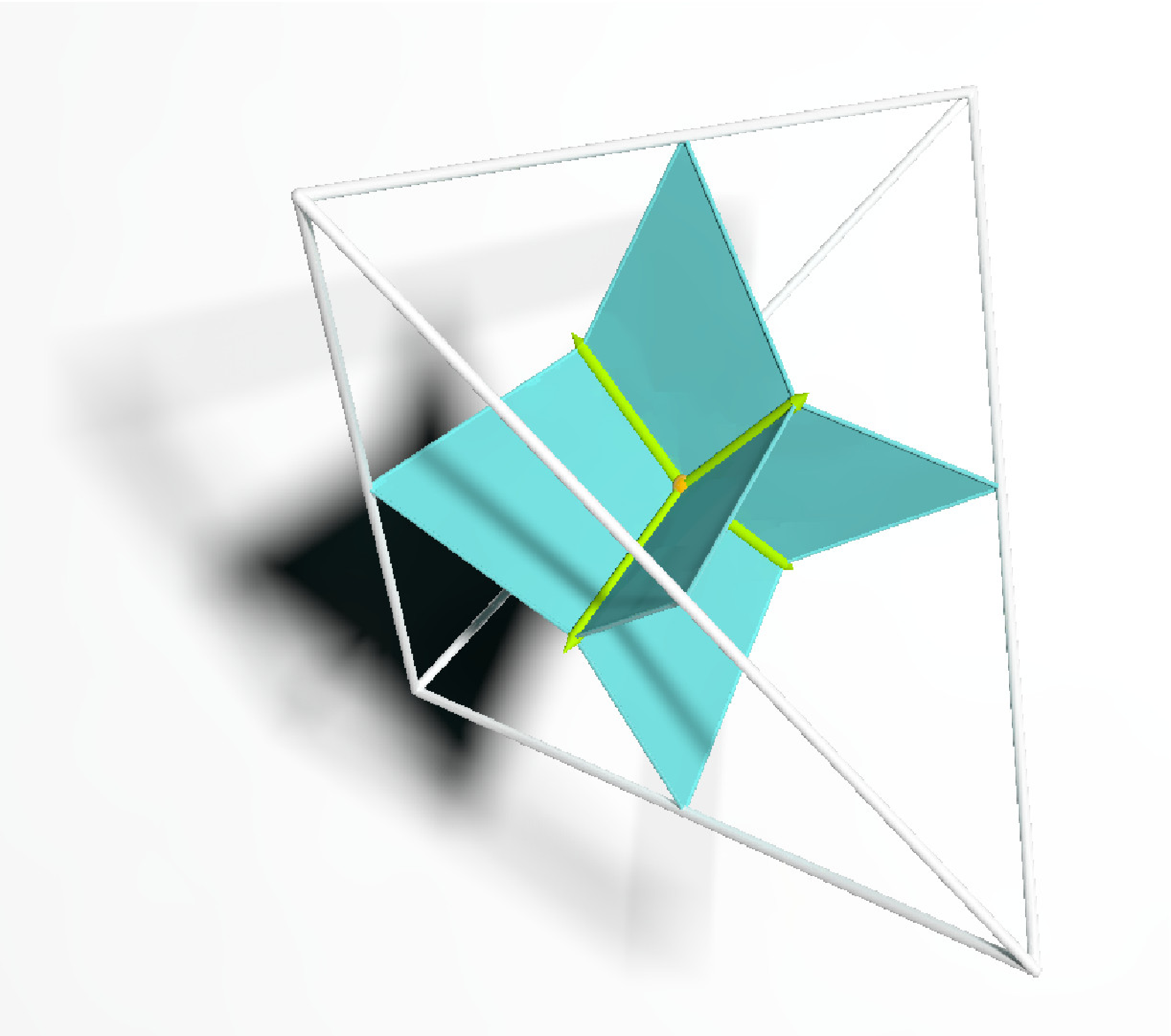}$$
\centerline{Dual polyhedral bits inside a tetrahedron $\Delta^3$}

Generally speaking, if $N$ is a triangulated PL-manifold, the tangent bundle $TN$ is
{\it not defined}, moreover, it is frequently not unique when it {\it is} defined \cite{HM}. 
Thankfully, non-smoothable PL structures, and distinct smoothings of PL-structures 
do not appear below dimension $7$. Thus the regular neighbourhoods of the dual $2$-skeleton
of a triangulated PL-manifold {\it do} have unique smoothings as the links of codimension
$1$ and $2$ faces are $0$-spheres and $1$-spheres respectively, which have unique
smooth structures, see \cite{HM}.  In particular, $TN_{|P^2}$ can be referred to without ambiguity and
we can discuss spin structures on PL-manifolds. 

\section{Geometry of simplices}\label{simpgeom}

This section describes some group-theoretic preliminaries related to the geometry of simplices.   
Let $\Sym(X)$ be the full group of isometries of an object $X$, and let $\Sym^+(X) \subset \Sym(X)$ be
 the orientation-preserving subgroup, provided these concepts make sense.    
Let $\mathcal D = \{(x_0,x_1,\cdots,x_n) \in \Real^{n+1} : x_0 = x_1 = \cdots = x_n \}$
be the `thin' diagonal, and let $\mathcal A = \{(x_0,x_1,\cdots,x_n) \in \Real^{n+1} : x_0+x_1+\cdots+x_n = 0\}$ 
be the anti-diagonal. 

Symmetries of $\Delta^n$ are determined by how they permute the vertices, thus there is an 
identification $\Sym(\Delta^n) \equiv \Sigma_{n+1}$ and $\Sym^+(\Delta^n) \equiv A_{n+1}$.  
If we translate $\Delta^n$ to the origin 
$$ \Delta^n_0 = \{ (x_0,\cdots,x_n) \in \mathcal A : x_i \geq \frac{-1}{n+1} \ \forall i \}$$
a linear extension gives an embedding $\Sym(\Delta^n_0) \to \gSO_{n+1}$.  The set $\mathcal D$ is an 
eigenspace relative to an eigenvalue $+1$ when the symmetry preserves orientation, and an eigenspace 
relative to an eigenvalue $-1$ when it reverses the orientation of $\Delta^n_0$ respectively. 

We will now examine the relations between the symmetric group and the {\it group of motions} of an
$n$-simplex.  Let $\Emb(\Delta^n, \Real^{n+1})$ be the space of affine-linear
embeddings of the $n$-simplex in $(n+1)$-dimensional Euclidean space.  
The space $\Emb(\Delta^n,\Real^{n+1})$ has the homotopy-type of a Stiefel manifold -- the displacement
vectors from one vertex to the remaining vertices gives such a map. This Stiefel manifold
in turn has the homotopy-type of $SO_{n+1}$ by Gram-Schmidt.
  
The group $\Sigma_{n+1}$ acts freely on the right on $\Emb(\Delta^n, \Real^{n+1})$ by re-labelling
the vertices of the simplex.  The group $\Sigma_{n+1}$ also acts on the left on $\Emb(\Delta^n,\Real^{n+1})$
by re-labelling the coordinate axes of $\Real^{n+1}$ but we will not need this action. 
The {\it motion group of the $n$-simplex} is defined to
be $\pi_1 \left( \Emb(\Delta^n, \Real^{n+1}) / \Sigma_{n+1} \right)$.  
Since $n \geq 2$ is always assumed, the homotopy long exact sequence of the bundle
$$\Sigma_{n+1} \to \Emb(\Delta^n, \Real^{n+1}) \to \Emb(\Delta^n, \Real^{n+1}) / \Sigma_{n+1}$$
gives us the $\Zed_2$-central extension
$$0 \to \Zed_2 \to \pi_1 \left( \Emb(\Delta^n, \Real^{n+1}) / \Sigma_n \right) \to \Sigma_{n+1} \to 0.$$

If $G$ is a group and $K$ an abelian group it is a standard theorem of group cohomology that the 
central extensions of $G$ with kernel $K$, taken up to extension-preserving isomorphism 
are in bijective correspondence with $H^2(G,K)$.  
It turns out that $H^2(A_n,\Zed_2)$ is a group of order two provided $n \geq 4$.  Thus, there is only 
one non-trivial $\Zed_2$-central extension of $A_n$.  
Schur called it the {\it double cover} of $A_n$, also called the {\it binary alternating group} and 
denoted either $2A_n$ or $\tilde A_n$. 
We use the latter notation. Schur also went on to show that $H^2(\Sigma_n,\Zed_2)$ is isomorphic to $\Zed_2^2$
for $n \geq 4$, moreover the restriction map $H^2(\Sigma_n,\Zed_2) \to H^2(A_n,\Zed_2)$ is onto,
thus there are two non-isomorphic $\Zed_2$-central extensions of 
$\Sigma_n$ which contain $\tilde A_n$.  We will give a geometric interpretation 
to one of these extensions.  A convenient notation for elements in these extensions
is given by Proposition \ref{schursthm}.

\begin{prop}\label{schursthm}\cite{Wilson} For all $n \geq 2$ there exist groups 
$\tilde \Sigma_n^+$ and $\tilde \Sigma_n^-$ which are $\Zed_2$-central extensions of $\Sigma_n$ 
such that:
\begin{enumerate}
\item Given a $k$-tuple $(a_1,\cdots,a_k)$ of distinct elements of $\{0,1,\cdots,n\}$ there is an 
element $[a_1 a_2 \cdots a_k] \in \tilde \Sigma_{n+1}^\pm$ called a $k$-cycle. 
\item The homomorphism $\tilde\Sigma_n^\pm \to \Sigma_n$ sends $[a_1 a_2 \cdots a_k]$ to 
$(a_1 a_2 \cdots a_k)$ for all $k$-cycles.
\item $[a_1 a_2 \cdots a_k] = [a_1 a_2 \cdots a_i][a_i a_{i+1} \cdots a_k]$ for all $k$ and
all $1 < i < k$.
\item If $\{a_1 a_2 \cdots a_k\}$ and $\{b_1 b_2 \cdots b_j\}$ are disjoint then
$[a_1 a_2 \cdots a_k][b_1 b_2 \cdots b_j] = $ 
\vskip 1mm $(-1)^{(k-1)(j-1)} [b_1 b_2 \cdots b_j][a_1 a_2 \cdots a_k].$
\item     $[a_1 a_2 \cdots a_k]^{[b_1 b_2 \cdots b_j]} =  
 (-1)^{(k-1)(j-1)}[\phi^{-1}(a_1) \phi^{-1}(a_2) \cdots \phi^{-1}(a_k)]$ where 
$\phi \in \Sigma_n$ is the cycle $(b_1 b_2 \cdots b_j)$.
We use the notation $g^h = h^{-1} g h$ for conjugation. 
\item  $[a_1 a_2 \cdots a_k]^k = \varepsilon$ for all $k \geq 2$, provided
$[a_1 a_2 \cdots a_k] \in \Sigma_{n+1}^\varepsilon$.
\end{enumerate}
\end{prop}

We call an element of $\tilde \Sigma_n^\pm$ {\it odd} or {\it even} if its projection to $\Sigma_n$ is odd or even respectively. 

Given that $\Sym(\Delta^n) \equiv \Sym(\Delta^n_0) \subset \gSO_{n+1}$, there is a canonical lift of 
$\Sym(\Delta^n)$ along the
2:1-covering map $\mSpin_{n+1} \to \gSO_{n+1}$.  We denote this 2:1-cover by 
$\widetilde{\Sym}(\Delta^n) \to \Sym(\Delta^n)$. It is a $\Zed_2$-central extension, since the
kernel of the map $\mSpin_{n+1} \to \gSO_{n+1}$ is central. 

\begin{prop}\label{prop1}
$\widetilde{\Sym}(\Delta^n)$ is canonically isomorphic to the motion group of $\Delta^n$ in $\Real^{n+1}$, 
$\pi_1 \left( \Emb(\Delta^n, \Real^{n+1}) / \Sigma_{n+1} \right)$.  It is also the $\Zed_2$ central extension
of $\Sigma_{n+1}$ denoted by $\tilde \Sigma_{n+1}^-$. Under this isomorphism, $\widetilde{\Sym^+}(\Delta^n)$ corresponds to $\tilde A_{n+1}$. 
\begin{proof} 
The isomorphism between the motion group $\pi_1 \left( \Emb(\Delta^n, \Real^{n+1}) / \Sigma_{n+1} \right)$ and the spin cover $\widetilde{\Sym}(\Delta^n)$ follows from the path-lifting property of the covering maps
$$\Sigma_{n+1} \to \Emb(\Delta^n, \Real^{n+1}) \to \Emb(\Delta^n, \Real^{n+1}) / \Sigma_{n+1}.$$  
Given an element of $\pi_1 \left( \Emb(\Delta^n, \Real^{n+1}) / \Sigma_{n+1} \right)$, lift a representative 
to a path in $\Emb(\Delta^n, \Real^{n+1})$ such that the endpoints differ by the action of $\Sigma_{n+1}$.  For such a lift, the initial embedding starts at the standard embedding of the simplex $\Delta^n$ in $\mathbb R^{n+1}$.  
Such a path extends to a path of affine linear automorphisms of $\mathbb R^{n+1}$, starting at 
${\mathrm {Id}}_{\Real^{n+1}}$. Using that $\gO_{n+1}$ is a deformation retract of $\GL(\mathbb R^{n+1})$ we can homotope this path (rel endpoints) to a path in $\gSO_{n+1}$, which therefore lifts to a path in $\mSpin_{n+1}$, starting at the identity element.  This describes the endpoint of the path as an element of $\widetilde{\Sym}(\Delta^n)$.

To verify that $\widetilde{\Sym}(\Delta^n)$ is isomorphic to $\tilde \Sigma_{n+1}^-$,
we will use the model $\Delta^n_0$ for the $n$-simplex. This 
has the advantage that the symmetries of $\Delta^n_0$ are linear. 
In this model, notice that either lift of the transposition 
$(a \ b)$ to $\mSpin_{n+1}$ has order $4$.  This is because in $\gSO_{n+1}$ the transposition 
$(a \ b)$ has a $2$-dimensional $(-1)$-eigenspace whose orthogonal complement is fixed pointwise. 
The $(-1)$-eigenspace is spanned by the vectors $e_a-e_b$ and $\sum_{i=0}^n e_i$. 
Thus if we denote any lift of $(a \ b)$ by $[a \ b]$ then $[a \ b]^2=-1$.  This is a 
proof by reduction to a universal example, as it is a direct computation to identify 
the spin cover $\mSpin_2 \to \gSO_2$ with the map of the unit
circle in the complex plane  $S^1 \to S^1$ given by $z \longmapsto z^2$. 
Relation 6 holds for $k=2$, and therefore for all $k \geq 2$. 
\end{proof}
\end{prop}

Consider the subgroup of $A_{n+1}$ which preserves the set $\{n-1, n\}$ i.e. its elements 
either fix $n-1$ and $n$ pointwise, or transpose them.  Since an element of $A_{n+1}$ 
is determined by its value on $n-1$ points, this subgroup is isomorphic to $\Sigma_{n-1}$. 
Thus, corresponding to a codimension-$2$ face of $\Delta^n$ there is an associated 
inclusion $\Sigma_{n-1} \to A_{n+1}$.   The lift of this $\Sigma_{n-1}$ to 
$\tilde A_{n+1}$ is isomorphic to $\tilde \Sigma_{n-1}^-$.  Using the notation of 
Proposition \ref{schursthm} one can verify that an embedding $\tilde \Sigma_{n-1}^- \to \tilde A_{n+1}$ is given by
$$ A \longmapsto \left\{ \begin{matrix} 
  A & A \text{ is even } \\
  A [n-1 \ n] & A \text{ is odd } \end{matrix} \right. .$$ 

There are precisely two embeddings $\tilde \Sigma_{n-1}^- \to \tilde A_{n+1}$ which cover 
the standard inclusion $\Sigma_{n-1} \to A_{n+1}$ ($\{0,1,2,\cdots,n-2\} \subset \{0,1,2,\cdots,n\})$.  
These two inclusions are essentially the same, as they differ by a
pre-composition with an automorphism of $\tilde \Sigma_{n-1}^-$ that fixes $\tilde A_{n-1}$ pointwise. 
The automorphism is given by
$$\tilde \Sigma^-_{n-1} \ni \sigma \longmapsto (-1)^{|\sigma|} \sigma \in \tilde\Sigma^-_{n-1}$$
 where $|\sigma|$ is the parity of $\sigma$. This is the unique non-trivial automorphism of
$\tilde \Sigma^-_{n-1}$ that fixes $\tilde A_{n-1}$ pointwise.

\section{Representing spin structures on triangulated manifolds}\label{mainsec}

As in Section \ref{defsnot}, let $P$ be the dual polyhedral decomposition to $T$, a 
triangulation of an $n$-manifold $N$.  We remind that the $k$-skeleton of $P$ is denoted $P^k$ while the 
set of $k$-cells is denoted $P_k$.  

This section gives a combinatorial technique to encode homotopy classes of sections of the Stiefel manifold of 
$(n-1)$-frames of $N$ over the dual $1$-skeleton $P^1$.  It also gives a combinatorial technique to determine which 
of these sections extend over the $2$-skeleton $P^2$. Once $N$ is oriented, this is our formalism 
for encoding spin structures as an $(n-1)$-frame extends to an oriented $n$-frame uniquely up to homotopy. 

To encode sections of the Stiefel bundle over $P^1$, we make the sections as
close to {\it simplicial} as possible.  That is:
\begin{itemize}

\item[1)] Evaluated at a point of $p \in P_0$ the vectors of each section should point towards
some of the vertices of the top-dimensional simplex containing $p$.    

\item[2)] We demand that the vectors of each section, evaluated at the barycentre of each dual edge $e \in P_1$ 
points to some of the vertices of the codimension one facet $F \in T_{n-1}$ dual to $e$. 
\end{itemize}

The total space of the Stiefel manifold of linearly-independent $(n-1)$-frames over $P^1$ will be
denoted by $V_{n-1} TN_{|P^1}$.  The set of sections of the bundle $V_{n-1} TN_{|P^1} \to P^1$ is written
$\Gamma_1 VN$, and the subspace satisfying (1) and (2) above is abbreviated by $\Sec$.

Every section in $\Gamma_1 VN$ is homotopic to one satisfying conditions (1) and (2).  Given any section, 
perform the homotopy along the finite set $P^0 \cup (T^{n-1} \cap P^1)$, and extend to $P^1$ via the 
homotopy extension property. 
To do this we use only that the Stiefel space $V_{n-1} \Real^n$ is connected, and that the 
`vertex pointing' subset of $V_{n-1} \Real^n$ is a non-empty set.  Thus the inclusion 
$\pi_0 \Sec \to \pi_0 \Gamma_1 VN$ is onto. 

$$\includegraphics[width=7cm]{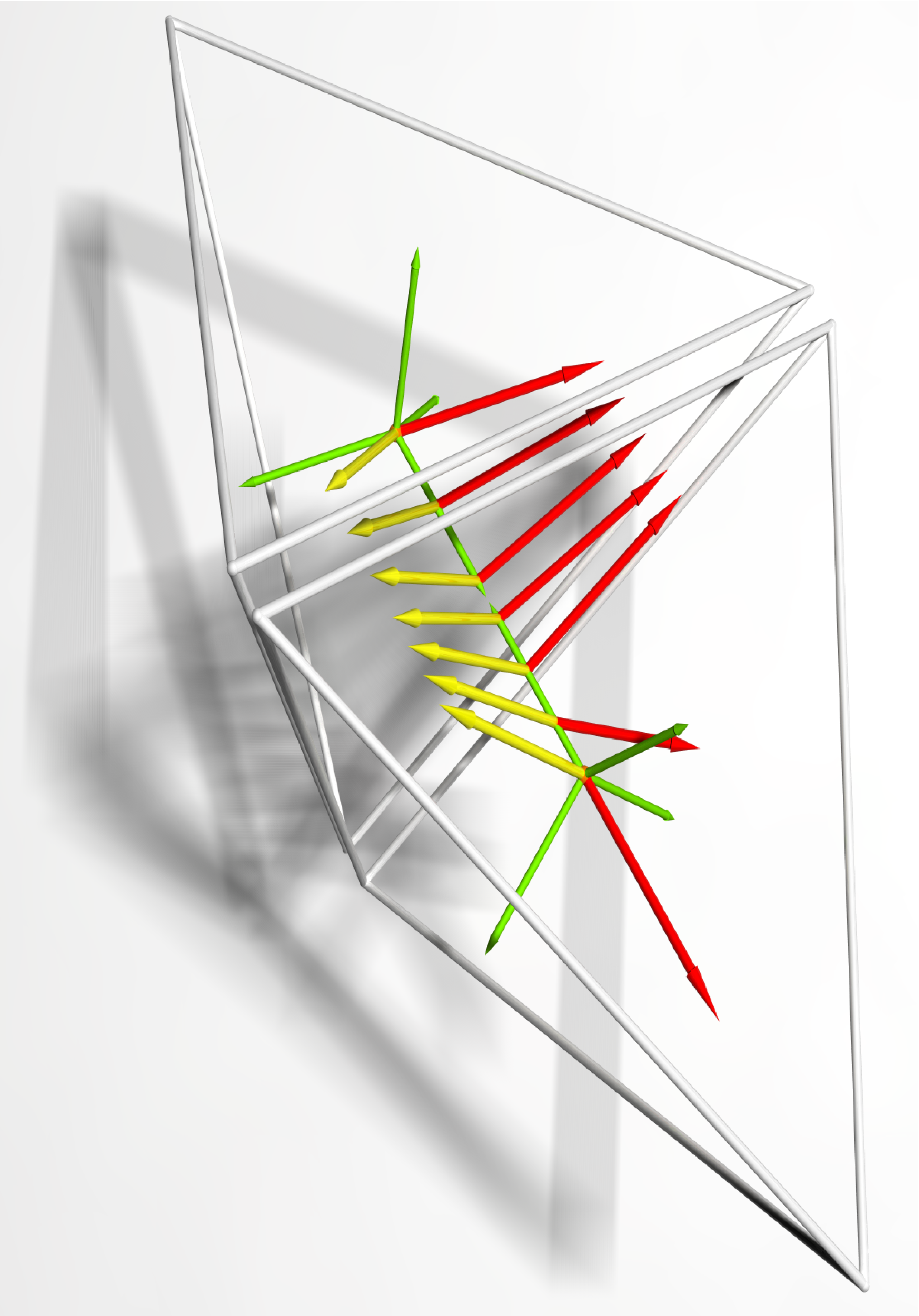}$$
\centerline{A section of $V_{n-1}(TN)_{|P^1}$ satisfying (1), (2).}
\centerline{The dual $1$-skeleton $P^1$ in green, $n=3$.}

We explain below how the set $\pi_0 \Sec$ can be thought of as a subset of the product
$$ \prod_{S \in T_n} A_{n+1} \times \prod_{F \in \sqcup_2 T_{n-1}} \tilde A_{n+1}. $$

\begin{itemize}
\item The $A_{n+1}$ corresponding to $S \in T_n$ factor encodes the section at the barycentre
of a top-dimensional simplex $S$. The
injections $\{0, 1, 2, \cdots, n-2\} \to \{0,1,2, \cdots, n\}$ have unique extensions 
to alternating bijections, so can be considered as elements of $A_{n+1}$.  
\item $\sqcup_2 T_{n-1}$ denotes two copies of $T_{n-1}$, one for each side of a dual edge $e \in P_1$ 
split at its barycentre.  Our conditions (1) and (2) above fix the behaviour of the section at the
end points of half of a dual edge $e$.  Thus corresponding to every dual edge there are 
two elements of $\tilde A_{n+1}$ that determine the element of $\pi_0 \Sec$.  We use the characteristic
map of $e \in P_1$ to determine which half of $e$ is the `first half' and which is the `second half'. 
\end{itemize}

The condition that sections are `face pointing' at the barycentre of a dual edge $e$ 
gives a constraint -- i.e. $\pi_0 \Sec$ is a {\it proper subset} of the above product.  We call this first
constraint the {\it face constraint}.  The continuity of our sections (at the barycentre of each edge $e \in P_1$) 
forces one further constraint, which we call the {\it continuity constraint.}

$$\pi_0 \Sec \subsetneq \prod_{S \in T_n} A_{n+1} \times \prod_{F \in \sqcup_2 T_{n-1}} \tilde A_{n+1}.$$

We call the set $\Sec$ the {\it simplicial sections} of the bundle $V_{n-1}(TN)_{|P^1}$.
We express the face and continuity constraints as formulas involving the characteristic 
maps of the triangulation, below. 

Given a codimension-$1$ face $F$ incident to a top-dimensional simplex $S$, let $\chi_F : \Delta^{n-1} \to N$
and $\chi_S : \Delta^n \to N$ be the characteristic maps respectively, and let $\iota \in \Sigma_{n+1}$ be the
{\it characteristic inclusion} of $F$ in $S$, i.e. $\chi_F = \chi_S \circ \iota$.  If $e$ is the edge dual to $F$ 
in $S$, let $\beta \in \tilde A_{n+1}$ be the motion corresponding to the half-edge $e \cap S$. 
The {\it face constraint} can be expressed as 
$$\beta(\alpha (\{0,1,\cdots,n-2\})) \subset \iota(\{0,1,\cdots,n-1\}),$$ 
where $\alpha \in A_{n+1}$ represents the section at the barycentre of $S$.  This is the formula
that says directly that the vertices pointed to by the vectors fields at $e \cap F$ are the vertices
of $F$.  Equivalently, we could say the vector field does not point to the vertex opposite to the
face, $\iota(n) \notin \beta(\alpha(\{0,1,\cdots,n-2\}))$.

Given a codimension-$1$ face $F$ incident to two top-dimensional simplices $S_1$ and $S_2$, 
let $\iota_1, \iota_2 \in \Sigma_{n+1}$ be 
the respective characteristic inclusions of $F$ in $S_1$ and $S_2$ respectively.  The {\it continuity constraint}
can be expressed as $(\iota_2^{-1} \beta_2 \iota_2)^{-1} \iota_1^{-1} \beta_1 \iota_1 \in \langle [n-1,n] \rangle$, 
where
$\langle [n-1,n]\rangle$ is the subgroup generated by $[n-1,n]$ in $\tilde A_{n+1}$.
This statement is equivalent to the statement that the two functions $\iota_k^{-1} \beta_k \iota_k$, as bijections of
the set $\{0,1,2, \cdots, n\}$ agree when restricted to the subset $\{0,1,2,\cdots,n-2\}$. The group 
$\Sigma_{n+1}$ acts naturally by conjugation on the group $\tilde \Sigma^-_{n+1}$, and the 
symbol $\iota_k^{-1} \beta_k \iota_k$ indicates the right action of $\iota_k$ on $\beta_k$.  

To encode the homotopy relation, we proceed by induction on the skeleton of $P^1$, i.e. first we perform the
homotopy on the $0$-skeleton $P^0$, and then we extend to $P^1$ using the homotopy extension property. Finally we perform
the homotopy on the edges of $P^1$, {\it leaving the end points fixed}.  The advantage of this perspective is
that it allows us to see that the homotopy relation as the orbit-space of a group action on
the simplicial sections $\pi_0 \Sec$.

The motions of an $n$-simplex $\Delta^n$ are given by $\tilde A_{n+1}$ (see Section \ref{simpgeom}). 
Given a simplicial section we represent it as an element
$(\prod_S \alpha_S, \prod_F \beta_F) \in \prod_{S \in T_n} A_{n+1} \times \prod_{F \in \sqcup_2 T_{n-1}} \tilde A_{n+1}$.  
Let $A \in \tilde A_{n+1}$ correspond to a motion of the $n$-simplex $S$, then 
the result of applying the motion $A$ to the simplicial sections at the barycentre of $S$, and extending to the entire
simplicial section is
$$A.\left(\prod_{S'} \alpha_{S'}, \prod_{F} \beta_F\right) = 
\left( \prod_{S'} \left\{ \begin{matrix} A\alpha_{S'} & \text{ if } S=S' \\
            \alpha_{S'} & \text{ if } S \neq S' \end{matrix} \right\}  , 
 \prod_F \left\{ \begin{matrix} \beta_F A^{-1} & \text{ if } F|S \\
           \beta_F & \text{ if } F \nmid S \end{matrix} \right\} \right), $$

where $F | S$ means `$F$ is a boundary facet of $S$' or equivalently, `$S$ is incident to $F$'.  This is
a group action of $\tilde A_{n+1}$ on $\pi_0 \Sec$. Moreover, observe that if $S_1$ and $S_2$
are distinct top-dimensional simplices of the triangulation $T$, then the two actions commute. 

To complete the description of the homotopy relation on simplicial sections, we describe the result of a homotopy 
of the section on the interior of an edge (fixed on the complement of the edge's interior). 
In principle, the justification for the formula below is the same as above, i.e. the standard algebra of 
obstruction theory, but the formula is made more complicated due to a change-of-coordinates issue. 
We have chosen to store all our motion data in the coordinates of the ambient top-dimensional simplex 
where the motion occurs.  But the edges of the dual cell complex $P_1$ {\it cross} from one top-dimensional
simplex to another, across a face $F \in T_{n-1}$. The group of motions of $F$ (fixing its barycentre) in the 
ambient triangulation is $\tilde \Sigma_n^-$, so we must provide the formalism for converting from the motions of
$F$ to motions in the adjacent top-dimensional simplices.  Given $A \in \tilde \Sigma^-_n$ representing a motion
of $F$ in the ambient triangulation (fixing the barycentre), the result of performing that homotopy on an element
of $\pi_0 \Sec$ at the barycentre of the edge, fixing the section outside the edge is given by

$$ A.\left( \prod_S \alpha_S, \prod_{F'} \beta_{F'}\right) = 
     \left( \prod_S \alpha_S, \prod_{F'} \left\{ 
\begin{matrix} \beta_{F'} & \text{ if } F' \neq F \\
      \iota A \iota^{-1} \beta_{F'} & \text{ if } F=F' \text{ and } $A$ \text{ even }  \\
      \iota [n \ A(k)]A \iota^{-1} \beta_{F'} & \text{ if } F=F' \text{ and } $A$ \text{ odd}
\end{matrix} \right\} \right).$$

In the above formula, $\iota$ is the characteristic inclusion of $F$ in $S$, and $k$ is the
index of the vertex in $F$ that is missed by the vector fields, i.e. 
$\{n, k\} = \iota^{-1} \beta_{F'} \alpha_S(\{n-1,n\})$.

There are a variety of ways to justify this formula, perhaps the most pragmatic is to consider the 
two $n$-dimensional simplices $S_1$ and $S_2$ incident to $F$ as two faces of some abstract $n+1$ dimensional
simplex $\mathcal S$ that is not part of the triangulation $T$ i.e. 
$S_i \subset \partial \mathcal S$ for $i = 1, 2$. If $\{0,1,2,\cdots,n-1\}$ is the vertex set for $F$, let 
$\{0,1,2,\cdots,n\}$ be the vertex set for $S_1$, $\{0,1,2,\cdots,n-1,n+1\}$ the vertex set for $S_2$
and $\{0,1,2,\cdots,n,n+1\}$ the vertex set for $\mathcal S$.  In Section \ref{simpgeom} we defined the inclusion
$\tilde \Sigma^-_n \to \tilde A_{n+2}$ via the formula $A \longmapsto [n+1 \ n]A$ provided $A$ is an odd
permulation. Notice that $[n+1 \ n]A$ maps $k \longmapsto A(k)$ and maps $n+1 \longmapsto n$.  So if we post-compose
$[n+1 \ n]A$ with the 3-cycle $[n \ n+1 \ A(k)]$ (which is the minimal motion in $\mathcal S$ returning $S_i$ to its initial position, or stated another way, this motion applied to $F$ projects into $N$ as an embedding) we get
$$[n \ n+1 \ A(k)][n+1 \ n]A = [A(k) \ n][n \ n+1]^2 A = [n \ A(k)]A.$$

We can replace $[n \ A(k)]$ with $-[n \ A(k)] = [A(k) \ n]$ in the above formula, as it simply corresponds to the
opposite embedding $\tilde \Sigma^-_{n-1} \to \tilde A_{n+1}$, which is just a convention for how lower-dimensional
motions convert to higher-dimensional motions.  

\begin{prop}\label{htpyprop} The maps defined above, for every $S \in T_n$
$$\tilde A_{n+1} \times \pi_0 \Sec \to \pi_0 \Sec$$
and for every $F \in T_{n-1}$
$$\tilde \Sigma_n^- \times \pi_0 \Sec \to \pi_0 \Sec$$
are group actions, moreover, the actions commute.  This gives us a group action
$$\left( \prod_S \tilde A_{n+1} \times \prod_F \tilde \Sigma_n^-\right) \times \pi_0 \Sec \to \pi_0 \Sec$$

whose stabilizers are isomorphic to $\Zed_2^m$ where $m$ is the number of path-components of $N$. 
The orbits of this action correspond to $\pi_0 \Gamma_1 VN$ via the map $\pi_0 \Sec \to \pi_0 \Gamma_1 VN$, 
thus if $N$ is oriented the orbits correspond canonically to spin structures. 
\begin{proof}

The maps for the top-dimensional simplices $\tilde A_{n+1} \times \Gamma^S_1 VN$ are group 
actions as they are essentially the canonical left and right actions of $\tilde A_{n+1}$ on $A_{n+1}$ and
$\tilde A_{n+1}$ respectively.

The maps involving $\tilde \Sigma^-_n$ corresponding to the codimension one faces $F$ require a more
subtle argument.  In the special case of even permutations, this is again the standard action of
$\tilde A_n$ on $\tilde A_{n+1}$, after conjugation by $\iota$.  

Let us consider the case where $A \in \tilde A_{n+1}$ can be an odd permulation. There are two non-trivial cases, 
$A_1.A_2.\beta$ where both $A_1$ and $A_2$ are odd, and the case $A_1$ odd and $A_2$ even.

(1) Let us first face the case where both $A_1$ and $A_2$ are odd.  If we let $k_2$ satisfy
$\{n, k_2\} = \iota^{-1} \beta \alpha \{n-1, n\}$ then we have
$$A_1.(A_2.\beta) = A_1.\left( \iota[n \ A_2k_2] A_2 \iota^{-1} \beta \right)$$
then $\{n, k_1\} = \iota^{-1} \iota[n \ A_2 k_2] A_2 \iota^{-1} \beta \alpha \{n-1 \ n\} = \{n \ A_2 k_2\}$, 
giving
$$A_1.A_2.\beta = \iota[n \ A_1 k_1] A_1 \iota^{-1} \iota [n \ A_2k_2] A_2 \iota^{-1} \beta$$
$$ = -\iota A_1 [n \ A_2 k_2][n \ A_2 k_2] A_2 \iota^{-1} \beta = \iota A_1 A_2 \iota^{-1} \beta = (A_1A_2).\beta.$$

(2) Now consider the case where $A_1$ is odd and $A_2$ is even.  The argument is simpler: 
$$A_1.(A_2.\beta) = A_1.(\iota A_2\iota^{-1} \beta)$$
$$ = \iota [n \ A_1 A_2 k]A_1 \iota^{-1} \iota A_2 \iota^{-1} \beta = \iota [n \ A_1A_2k]A_1A_2\beta = (A_1A_2).\beta.$$

We now establish that the kernel is $\Zed_2^m$.  If an element is stabilized 
under the action, the underlying section at $P^0$ and motions along $P^1$ 
are fixed.  This forces the components of the $\prod_S \tilde A_{n+1}$ factor to be $\pm 1$. 
Similarly the components of the product $\prod_F \tilde \Sigma^-_n$ are all $\pm 1$.  Thus
we can think of the elements of the stabilizers as $0$-dimensional mod-2 cocycles.  
Such objects correspond to $H^0(N, \Zed_2)$ which is isomorphic to $\Zed_2^m$. 
\end{proof}
\end{prop}

We turn our attention to spin structures -- the issue of determining which elements of
$\Gamma^S_1 VN$ admit extensions to $P^2$.  Let $W \in T_{n-2}$ be a codimension $2$ simplex of $T$. 
Let $S_0, S_1 \cdots, S_{m-1}$ be the circuit of $n$-simplices about $W$.  
Let $F_0, F_1, \cdots, F_{m-1}$ be the corresponding circuit of codimension $1$ simplices.   
We choose these coherently -- the normal sphere to $W$ in $T$ is a triangulated circle and we 
index the $S_i$ and $F_i$'s in accord with that cyclic order. 

Given an element of $\pi_0 \Sec$ we will set up a formula representing
the obstruction to extending it over the $2$-cell dual to $W$.  
Let $\beta_{1 i}$ and $\beta_{2 i} \in \tilde A_{n+1}$ be the motions of
the simplex $S_i$ as one travels from the barycentre of $S_i$ to the barycentres of the faces 
$F_{i-1}$ and $F_i$ respectively (the index $i$ taken mod $m$).  
Our perspective will be to {\it cut} the normal circle to $W$ in $T$ and align the simplices $S_0, S_1, \cdots, S_{m-1}$
as if they were parallel.  We compute the motion of the vector fields as one traverses the circuit of simplices, 
in these parallelized coordinates. Let $w_i : \Delta^{n-2} \to \Delta^n$ be the characteristic inclusion corresponding
to $W \hookrightarrow S_i$, extending uniquely to be an element $w_i \in A_{n+1}$ via the condition that
$w_i(n-1)$ and $w_i(n)$ represent the vertices of an edge of the normal circle to $W$, with its cyclic orientation. 
The {\it parallelized total motion} in the simplex $S_i$ (about $W$) will be denoted by $S_i^w$, is defined as

$$S_i^w = \left\{ \begin{matrix} 
 \beta_{2i}^{w_i} (\beta_{1i}^{w_i})^{-1} & 
   \text{ if } [0, n-2] \setminus w_i^{-1} \beta_{1i} \alpha_i = \emptyset \\
 [a \ n \ n-1] \beta_{2i}^{w_i} (\beta_{1i}^{w_i})^{-1} & 
   \hskip 4.4mm \text{ if } [0, n-2] \setminus w_i^{-1} \beta_{1i} \alpha_i = \{a\}
\end{matrix}\right. .$$

\begin{prop}\label{extthm}
An element of $\pi_0 \Sec$ extends over the $2$-cell dual to $W \in T_{n-2}$ 
if and only if the product of the parallelized total motions is the non-trivial central element of $\tilde A_{n+1}$, 
namely

$$S_{m-1}^wS_{m-2}^w\cdots S_1^wS_0^w = -1.$$
\end{prop}

The explanation for this formula is in a similar spirit to Proposition \ref{htpyprop}. Imagine two consecutive 
simplices $S_i$ and $S_{i+1}$ stuck together along their common face $F_i$, and imagine the simplices 
pulled-apart so that they are parallel.  The motion $\beta_{2i}^{w_i} (\beta_{1i}^{w_i})^{-1}$ is what one
applies to the vector fields in $S_i$ as one travels along $P^1$ from the face $F_{i-1}$ to $F_i$ in the 
simplex $S_i$.  Consider what additional motion we need to apply to these vector fields as we rotate 
$S_{i+1}$ to be parallel to $S_i$. If our vector fields at $F_i$ are pointing into $W$, we would be 
done because the motion that makes the simplex $S_{i+1}$ parallel to 
$S_i$ has no affect on the vector fields. This is the case when 
$[0, n-2] \setminus w_i^{-1} \beta_{1i}\alpha_i = \emptyset$.  
If the vector fields hit the vertex of $F_i$ not in $W$, i.e. they miss a vertex of $W$, then
$[0, n-2] \setminus w_i^{-1}\beta_{1i}\alpha_i = \{a\}$ and our motion to rotate $S_{i+1}$ to be parallel with $S_i$ affects the vector fields.  The motion can be expressed as $[a \ n \ n-1]$ in the coordinates of $F_i$, hence the formulas for $S_i^w$. Our formulas for $S_i^w$ with $i=0, 1, 2, \cdots, m-1$ are now in a common `parallel' coordinate system and can be concatenated. 
We demand the product is $-1$ since the act of `closing' the parallel simplices contributes an extra $2 \pi$ rotation into the product. 

Thus, a {\it combinatorial spin structure} on a triangulated, oriented $n$-manifold $N$ is an orbit of
$\prod_S \tilde A_{n+1} \times \prod_F \tilde \Sigma_n^-$ acting on 
$\pi_0 \Sec$, whose elements extend over all dual $2$-cells $W \in P_2$. 

\begin{eg}\label{eg1} In the diagram below the red arrows indicate the vector field
over the $0$-skeleton, given by $\alpha_i$, as well as the vector field when pushed into the faces $F_i$. 
The blue arrows indicate our convention that our motions are specified as
motions as one travels {\it from} the barycentres of top-dimensional simplices {\it to} the 
barycentres of the codimension one simplices $F_i$.   We have chosien to embed the triangles in the
plane so that $(012)$ represents a counter-clockwise $2\pi/3$ rotation.
\begin{center}
\begin{multicols}{2}
{
\psfrag{S1}[tl][tl][0.9][0]{$S_0$}
\psfrag{S2}[tl][tl][0.9][0]{$S_1$}
\psfrag{S3}[tl][tl][0.9][0]{$S_2$}
\psfrag{S4}[tl][tl][0.9][0]{$S_3$}
\psfrag{S5}[tl][tl][0.9][0]{$S_4$}
\psfrag{F1}[tl][tl][0.9][0]{$F_0$}
\psfrag{F2}[tl][tl][0.9][0]{$F_1$}
\psfrag{F3}[tl][tl][0.9][0]{$F_2$}
\psfrag{F4}[tl][tl][0.9][0]{$F_3$}
\psfrag{F5}[tl][tl][0.9][0]{$F_4$}
\psfrag{W}[tl][tl][0.9][0]{$W$}
$\includegraphics[width=6cm]{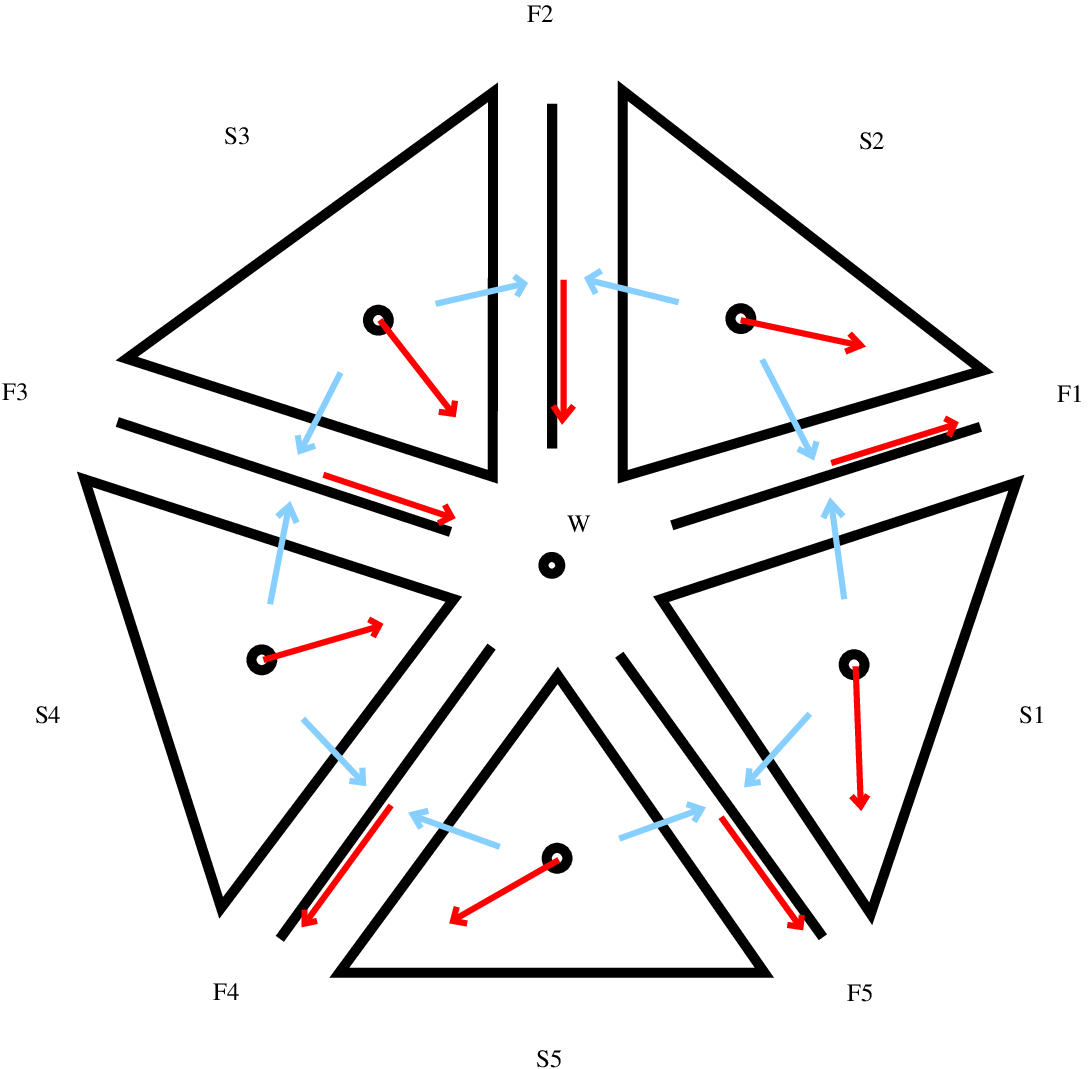}$
}
\columnbreak
\vskip 0.5cm
$\begin{aligned}
& S_0^w = \pm 1 \\ 
& S_1^w = \pm [210]\\ 
& S_2^w = \pm 1,   \\ 
& S_3^w = \pm [012]\\ 
& S_4^w = \pm 1.
\end{aligned}$ 
\end{multicols}
\end{center}
\end{eg}

In Example \ref{eg1}, the if we choose the `short' motions consistent with the figure, i.e. all non-trivial turns are 
either clockwise or counter-clockwise by $\frac{2\pi}{3}$.  We complete the computation with 
$$S_0^w = 1, \ S_1^w = [021], \ S_2^w = 1, \ S_3^w = -[012], \ S_4^w = 1,$$ giving
$$S_4^w S_3^w S_2^w S_1^w S_0^w = -[012][021] = -1.$$
By Proposition \ref{extthm} (or by visual inspection) the vector field extends over the $2$-cell dual to $W$.

Although the parameter space for $\pi_0 \Sec$ has order
$$ |A_{n+1}|^{|T_n|}|\Sigma_n|^{|T_{n-1}|}4^{|T_{n-1}|} $$
one could implement this formalism by assuming the vector fields over $P^0$ are induced by the characteristic maps, similarly for the vector fields on the barycentres of the dual edges $e \in P_1$.  One can assume that on one half of each edge the vector fields chosen are given by some canonical path. In this setup $\pi_0 \Sec$ is a subset of a set parametrized by $|F|$ bits.  
This is analogous to orientations: orientation can be thought of as a plus or minus sign $\pm$ associated to every top-dimensional simplex, satisfying a coherence condition. Spin structures are similarly parametrized by a $\pm$ sign on every dual edge $e \in P_1$ satisfying an analogous coherence condition. 

The first Stiefel-Whitney class $\omega_1 \in H^1(N,\Zed_2)$ of a manifold $N$ is the obstruction to orientability.  From the perspective of triangulations, the $1$ cocycle representing $\omega_1$ is given by comparing the orientations of top-dimensional simplices adjacent across a face $F$.  If they are oriented compatibly, meaning the transition function $\phi$ satisfies $\phi \in \Sigma_{n+1} \setminus A_{n+1}$, then $\omega_1(F) = 1$, otherwise $\omega_1(F)=-1$.  

There is a similar computation of $\omega_2$, the second Stiefel-Whitney class.  As a 2-cocycle, $w_2$ is computed by constructing 
$n-1$ everywhere linearly independent sections on $P^1$.  Its value on a $2$-cell dual to $W \in T_{n-2}$ is precisely 
our extension obstruction $-S_{m-1}^wS_{m-2}^w\cdots S_1^wS_0^w$.

\section{Combinatorial complex spin structures}

As described in Section \ref{defsnot}, a \tspinc-structure on an $n$-manifold $N$ consists of a homotopy class of a lift of the tangent bundle classifying map $N \to B\gO_n$ to the $\mSpinc$ classifying space $B\mSpinc_n$.  We take
the perspective of Section \ref{defsnot} and consider a \tspinc-structure on $N$ as a $1$-dimensional complex bundle over $N$, call it $\nu$, together with a spin structure on the sum of the two bundles $TN \oplus \nu$.  One-dimensional complex bundles over $N$ are classified by maps $N \to B\gSO_2$, which correspond precisely (via obstruction theory) to elements of $H^2(N,\Zed)$. 

Let $\beta$ be a co-chain representing an element of $H^2(N,\Zed)$, $\nu$ the complex line bundle associated to $\beta$, and consider the problem of finding a spin structure on $TN \oplus \nu$.  Since $\pi_1 \gSO_3 \to \pi_1 \gSO_5$ is an isomorphism, we can demand that our trivialization of $TN \oplus \nu$ over $P^1$ is the direct sum of a trivialization of $TN$ over $P^1$ with a fixed trivialization of $\nu$ over $P^1$. Since $\gSO_2$ is connected $\nu$ is trivial when restricted to $P^1$. 

Checking whether or not such a trivialization of $TN \oplus \nu$ over $P^1$ extends to a trivialization over $P^2$, we get the condition
$$S_{m-1}^wS_{m-2}^w\cdots S_1^wS_0^w = (-1)^{1+ \beta(W)}$$
where $\beta(W)$ is the value of $\beta$ on the $2$-cell dual to $W$, and the remainder of the formula is as
in Section \ref{mainsec}.

Thus by design, $N$ has a combinatorial \tspinc-structure if and only if $\omega_2$ is the mod-2 reduction of a class in $H^2(N,\Zed)$. More specifically, there exists a spin structure on $TN \oplus \nu$  if and only if $\omega_2$ is the mod two reduction of $\beta$. 

\section{Appendix}

This section collects a few lesser-known facts related to the paper. These results are useful to anyone
interested in implementing these techniques in software, and they are available in the software package Regina \cite{BBP}.

When $n$ is odd, $\Sym(\Delta^n)$ has an alternative interpretation. 
There is a canonical isomorphism $\Sym(\Delta^n_0) \simeq \Sym^+(\Delta^n_0 \cup -\Delta^n_0)$
where $-\Delta^n_0$ is the antipodal simplex.  The isomorphism is given by
$$\Sym(\Delta^n_0) \ni A \longmapsto (-1)^{|A|} A \in \Sym^+(\Delta^n_0 \cup -\Delta^n_0).$$ 

Given a subgroup $G$ of $\gSO_n$, let $\tilde G \subset \mSpin_n$ be the preimage of 
$G$ under the covering map $\mSpin_n \to \gSO_n$.  By design $\tilde G$ is a 
$\Zed_2$-central extension of the group $G$.

$$\xymatrix{0 \ar[r] & \Zed_2 \ar[r]        & \mSpin_n \ar[r]         & \gSO_n  \ar[r] & 0\\
	    0 \ar[r] & \Zed_2 \ar[r] \ar[u]^{\simeq} & \tilde G \ar[r] \ar@{->}[u] & G \ar@{->}[u] \ar[r] & 0}$$

\begin{center}
\scalebox{0.24}{\includegraphics{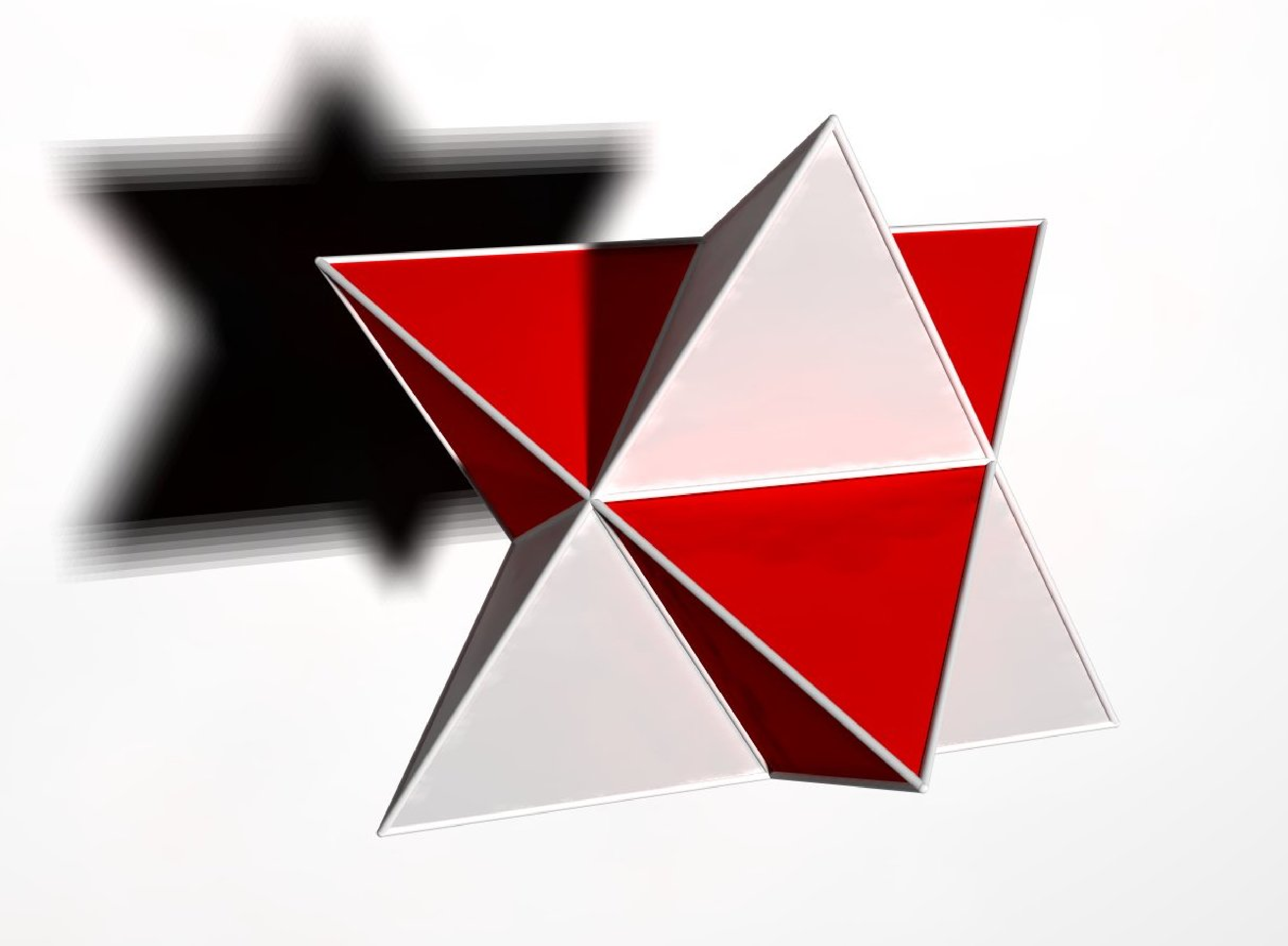}}
$$\Delta^3_0 \cup -\Delta^3_0$$
\end{center}

We give some concrete descriptions of the low-dimensional groups $\tilde A_5$ and
$\tilde \Sigma_4^-$ respectively. Although not required for the main results of the
paper, we compile this information here for easy reference. 

\begin{prop}\label{prop2}
The group $\mSpin_3$ has a natural identification with the unit sphere in the quaternions,
$S^3$. $S^3$ acts on the quaternions by conjugation.  This action is an orthogonal linear map, and
it fixes the real line pointwise.  If we call the orthogonal complement of $\Real$ in the quaternions
the {\it purely imaginary} quaternions, we can identify the purely imaginary quaternions with $\Real^3$. 
Thus our action can be interpreted as a Lie group homomorphism $S^3 \to \gSO_3$. 

Consider $\Delta^3$ to be the convex hull of the four points 
$$\left\{ 
     \left(\frac{-1}{\sqrt{3}}, \frac{1}{\sqrt{3}}, \frac{1}{\sqrt{3}}\right), 
     \left(\frac{1}{\sqrt{3}},  \frac{-1}{\sqrt{3}},  \frac{1}{\sqrt{3}}\right), 
     \left(\frac{1}{\sqrt{3}}, \frac{1}{\sqrt{3}},  \frac{-1}{\sqrt{3}}\right), 
     \left(\frac{-1}{\sqrt{3}},  \frac{-1}{\sqrt{3}}, \frac{-1}{\sqrt{3}}\right) \right\} \subset \Real^3.$$
Then $\tilde \Sigma_4^-$ is isomorphic to the subgroup of $S^3$ which preserves $\Delta^3 \cup -\Delta^3$.
It consists of the elements:
$$\left\{ \pm 1, \pm a, \frac{1}{\sqrt{2}}(\pm 1 \pm a), \frac{1}{\sqrt{2}}(\pm a \pm b), 
	\frac{1}{2}(\pm 1 \pm a \pm b \pm c) \text{ where } \{a,b,c\} = \{i,j,k\} \right\} \subset S^3.$$

The group $\mSpin_4$ has a natural identification with $S^3 \times S^3$. The
homomorphism $S^3 \times S^3 \to \gSO_4$ given by left and right multiplication by unit quaternions.

Consider $\Delta^4$ to be the convex hull of the points
$$\left\{ (1,0,0,0), 
\frac{1}{4}(-1, -\sqrt{5}, \sqrt{5}, \sqrt{5}), 
\frac{1}{4}(-1, \sqrt{5}, -\sqrt{5}, \sqrt{5}),
\frac{1}{4}(-1, \sqrt{5}, \sqrt{5}, -\sqrt{5}),
\frac{1}{4}(-1, -\sqrt{5}, -\sqrt{5}, -\sqrt{5}) \right\}.$$
The $120$ elements of $\tilde A_5 \subset S^3 \times S^3$ are given by: $\pm (1,1) $ 
having orders $1$ and $2$ respectively, $\left( \alpha, \alpha\right)$
where $\alpha = \frac{\pm 1 \pm i \pm j \pm k}{2}$, having orders $3$ and $6$ respectively, together 
with $\pm \left( \frac{1}{2} + a\alpha + b\beta, \frac{1}{2} + \overline{a}\alpha + \overline{b}\beta \right)$
where $a = \pm \frac{\sqrt{5}-1}{4}$, $b = \pm \frac{\sqrt{5}+1}{4}$, $\{\alpha, \beta, \gamma\} = \{i,j,k\}$,
and $\alpha\beta = \gamma$, where $\overline a$ indicates the image of $a$ under the automorphism of
$\Rat[\sqrt{5}]$ given by $\sqrt{5} \longmapsto -\sqrt{5}$. These elements have order $3$ and $6$ respectively.
There are also the elements 
$\pm \left( a \pm \overline{a}\alpha \pm \frac{\beta}{2}, \overline{a} \pm a\alpha \pm \frac{\beta}{2}\right) $
where $\{\alpha,\beta,\gamma\} = \{i,j,k\}$.  If $\alpha \beta = \gamma$ $a =  \frac{1-\sqrt{5}}{4}$, 
otherwise $a =  \frac{1+\sqrt{5}}{4}$.  These elements have order $5$ and $10$. 
There are the elements $ \pm (\alpha, \alpha) $ where $\alpha \in \{i,j,k\}$. These elements have order $4$.
Finally there are the elements $ \left( a\alpha + b\beta + c\gamma, 
		 \overline{a}\alpha+\overline{b}\beta+\overline{c}\gamma\right)$
where $\{\alpha,\beta,\gamma\} = \{i,j,k\}$, $\alpha\beta=\gamma$, $a = \pm \frac{1+\sqrt{5}}{4}$, 
$b=\pm \frac{1-\sqrt{5}}{4}$, $c = \pm \frac{1}{2}$.  These elements have order $4$.
\end{prop}

\providecommand{\bysame}{\leavevmode\hbox to3em{\hrulefill}\thinspace}

\Addresses

\end{document}